\documentclass[11pt]{article}
\usepackage{graphicx}        
\usepackage{amssymb}
\usepackage[numbers]{natbib}
\usepackage{bm}        
\usepackage{tikz}
\usepackage{smartdiagram}

\def\x{{\boldsymbol{x}}}
\newcommand{\brak}[1]{\left\{\begin{array}{llllllll} #1 \end{array}\right. }

\def\be{\begin{equation}}
\def\ee{\end{equation}}

\begin{document}

\title{Ten New Benchmarks for Optimization}
\author{Xin-She Yang \\ 
School of Science and Technology, \\ Middlesex University London, \\
The Burroughs, London NW4 4BT, United Kingdom.  }

\date{}

\maketitle

\abstract{Benchmarks are used for testing new optimization algorithms and their variants to evaluate their performance. Most existing benchmarks are smooth functions. This chapter introduces ten new benchmarks with different properties, including noise, discontinuity, parameter estimation and unknown paths. }

\bigskip 


\bigskip 
\noindent {\bf Citation Details:}
Xin-She Yang, Ten New Benchmarks for Optimization, in: {\it Benchmarks and Hybrid Algorithms in Optimization and Applications} (Edited by Xin-She Yang), Springer Tracts in Nature-Inspired Computing, pp. 19 -- 32 (2023). \\ 
https://doi.org/10.1007/978-981-99-3970-1\_2
\\[15pt]

\section{Introduction}

The literature of nature-inspired algorithms and swarm intelligence is expanding rapidly, and most nature-inspired optimization algorithms are swarm intelligence based algorithms~\citep{YangBook2014,YangHe2019,Yang2018Wiley,Yang2020Rev}.

New algorithms for optimization appear regularly in the current literature and hundreds of papers are being published every year, about new algorithms and their variants as well as their applications. Obviously, new algorithms have to be tested and validated using various known problems with known optimality locations. Simple benchmark functions are often the first set of problems that new algorithms or variants attempt to solve.

However, most existing benchmarks are function benchmarks that are usually smooth with known optimality at a single location for a given function benchmark. In addition, the search domains of these functions in terms of their independent variables are usually regular, often expressed in terms of simple bounds and limits. Even though some function benchmarks have constraints, these constraints tend to be sufficiently simple, which do not change the shape of the search domains significantly. The tests using such smooth benchmarks with regular search domains may give some insights into the performance of the algorithms under consideration. However, the usefulness of such benchmarks may be quite limited because these functions have almost nothing to do with the realistic optimization problems from the real-world applications.

Ideally, we should have a diverse set of test benchmarks and case studies derived from real applications so that researchers can use them for testing algorithms, especially new algorithms. However, such test benchmarks are largely not available because case studies tend to be very specialized in a specific subject area. Even we may have such case studies as benchmarks, specialized knowledge in a subject area is needed to solve the problem properly and interpret the results correctly. If we can somehow extract the essential part of the optimization problems and try to make them almost independent of special subjects, we can make the relevant problems as generic benchmarks.

In the rest of the chapter, we will first briefly discuss the role of benchmarking and the different types of benchmarks. Then, we will introduce ten new benchmarks with additional properties, such as noise, discontinuity, non-differentiability, multi-layered discrete values, and optimal paths from the calculus of variations.

\section{Role of Benchmarks}

There are many benchmark functions that are smooth functions with known optimal solutions and optimal objective values~\citep{Yang2010Fun,Jamil2013,Suganthan2005}. Such benchmarks enable researchers to test new algorithms so as to gain better understanding about the convergence behaviour, stability and performance of the algorithm under consideration.

Benchmarks can be very diverse. To validate any new optimization algorithm, a variety of test benchmarks should be used to see how the algorithm under consideration may perform for different types of problems. In general, benchmarks can be divided into five categories:

\begin{figure}[h]
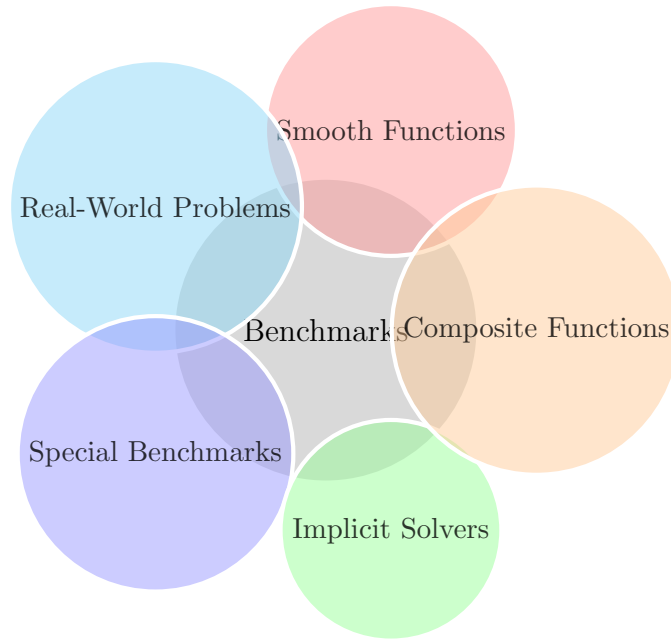

\begin{center}
\smartdiagramset{planet text width=2.5cm, satellite text width=2cm}
\smartdiagram[bubble diagram]{Benchmarks, Smooth Functions, Real-World Problems, Special Benchmarks, Implicit Solvers, Composite Functions}
\caption{Different types of benchmarks.}
\end{center}
\end{figure}

\begin{itemize}

\item \emph{Smooth Functions}: There are more than 200 smooth functions as benchmarks in the current literature~\citep{Jamil2013,MazharFun,HedarFun,Suganthan2005,Kumar2020}.
    Whether the problems are constrained or often unconstrained, their objective landscapes tend to be smooth. In addition, in some rare cases, such test functions can also be multi-objective optimization~\citep{Zitzler2000}.

\item \emph{Composite Functions}: Though many simple test functions such as the Ackley function exist, composite functions have been designed to make it harder for algorithms to find the optimal solutions because the locations of the optimality of these functions are shifted and twisted.

\item \emph{Implicit Solvers}: In many applications, the exact form of the objective can be difficult to express explicitly. For example, many engineering design problems require to use finite element analysis (FEM) packages and computational fluid dynamics (CFD) packages to evaluate the design performance. In this case, the actual evaluation of the design objectives are carried out by calling external solvers for given inputs (design variables and parameters) and then extracting the output (some objective or performance metrics).

\item \emph{Special Benchmarks}: Sometimes, specialized benchmarks are used to test certain types of optimization techniques. Benchmarks can be very specialized. For example, in protein folding and molecular biological applications, the test cases are often very specialized with given structure of data and objectives.

\item \emph{Real-World Problems}: For new optimization algorithms to be truly useful, they should be able to solve real-world problems. Therefore, extensive tests should be carried out using real-world problems as benchmarks. This category of benchmarks can be extremely diverse to cover almost all areas of applications.

\end{itemize}

Since there many different types of benchmarks, an important question is naturally: What benchmarks should be used for validating new algorithms?

For most of the studies in the current literature, the benchmarking process seems to use a finite set (typically about a dozen to two or three dozens of functions for testing algorithms, even though these functions are selected to have some diverse properties, such as mode shapes, optimality locations, separability of different dimensions, and even with some constraints. Though these functions themselves can be quite complicated with multiple modes, they are idealized functions, which may have nothing to do with problems arisen from real-world applications. Thus, whatever the conclusions may be drawn from testing functions, they may not be much use in practice. After all, most real-world problems can be even more complex with highly nonlinear constraints and irregular search domains. Therefore, it can be expected that algorithms work for test functions may not work well in practical applications.

\section{New Benchmark Functions}

Even the usefulness of simple benchmarking is quite limited, it is still an important part of the algorithm evaluation process. In addition to the existing test functions, we now introduce even more complex functions as new benchmarks. These new benchmarks highlight the challenges in algorithm testing and may inspire more research to test new algorithms from a wider perspective, including introducing noise, non-unique optimal solutions and even solutions in an infinite-dimensional functional space.

\subsection{Noisy Functions}

As almost all existing benchmark problems are deterministic in the sense that the function forms and their solutions have no randomness.
Now let us add some noise to a smooth function, which may make it more challenging for algorithms to find its optimality.

One way for adding noise without affecting the location of its optimality is to multiplying a random variable drawn from a uniform distribution. For example, we can have a noisy function
\be f(\x)=\sum_{n=1}^D \epsilon_n x_n^2, \ee
as the extension to a standard sphere function. Here, all $\epsilon_n$ are drawn from a uniform distribution in [0,1].

We can design an even more complicated function
\be f_1(\x)=\sum_{n=1}^D x_n^{2 n} -\exp\Big[-\sum_{n=1}^D \epsilon_n x_n^{2n} \Big], \ee
with \be  -100 \le x_n \le 100, \ee
where $D \ge 1$ is the dimensionality of the function. Again, all $\epsilon_n$ are drawn from a uniform distribution in [0, 1].
Its optimality $f_{\min}=-1$ occurs at $\x=(0,0, ..., 0)$.

\subsection{Non-differentiable Functions}

The simplest function with a kink is probably $f(x)=|x|$, which has a global minimum $f_{\min}=0$ at $x_*=0$. However, this function does not have a well-defined derivative at $x=0$. In the $D$-dimensional space, we can extend it to
\be f(\x)=\sum_{n=1}^D |x_n|, \ee
and its global minimum $f_{\min}=0$
is located at $\x=(0, 0, ..., 0)$.

Obviously, we can design more functions with multiple kinks, such as
\be f_2(\x) =\left(\sum_{n=1}^D |x_n-n \pi|\right) \exp\left[-\sum_{n=1}^D \Big|\sin (|x_n-n \pi|)\Big|  \right], \ee
with \be -D \pi  \le x_n \le +D \pi. \ee
Its optimality $f_{\min}=0$ occurs at $\x=(\pi, 2 \pi, ..., n \pi)$.

\subsection{Functions with Isolated Domains}

In most benchmark problems, their search domains are typically regular in the form $x_n \in [a,b]$, ranging in an interval from $x_n=a$ to $x_n=b$. This is usually true for unconstrained optimization problems with simple bounds or limits. For constrained optimization problems, depending on the actual constraints, feasible search domains can have very irregular shapes or even with isolated, fragmental regions.

For example, we can design a function with two isolated domains
\be f_3(\x)=x_1^2 + \sum_{n=2}^D |x_n^3|, \ee
subject to
\be |x_1-2 a|+\sum_{n=2}^D |x_n| \le a, \ee
and
\be \sum_{n=1}^D (x_n-5a)^2 \le a^2, \quad a \ge 1. \ee
The minimum $f_{\min}=a^2$ occurs at $\x=(a,0,0,...,0)$.

More generally, we can have a function with multiple isolated domains with four peaks as its optimal solutions~\citep{Yang2020Rev}
\be f_4(\x)=\sum_{i=-N}^N \;\; \sum_{j=-N}^N \Big(|i|+|j|\Big)
\exp \left[-a(x-i)^2-a(y-j)^2 \right], \ee
in the domains of
\be |x-i|+|y-j| \le \frac{1}{a}, \quad \textrm{for } -N \le i,j \le N, \ee
where $i,j$ are integers with $N=100$ and $a=10$. This function has $4(N+1)^2$ local peaks, but it has four highest peaks at four corners. However, its domain is formed by many isolated regions, or
$4(N+1)^2=40 401$ regions when $N=100$.

\section{Benchmarks with Multiple Optimal Solutions}

Almost all benchmark functions have their optimal objective values at a finite number of isolated points. For example, $f(x,y)=x^2+y^2$ has only a single optimal solution at the point $(0,0)$ with $f_{\min}=0$. To make things more complicated, we can design functions with infinitely many solutions with equal objective values. For example, the function $g(x,y)=x^2+y^2-2xy=(x-y)^2$ has the global minimum $f_{\min}=0$ on the line $y=x$.

Many researchers use the closeness to the optimal solution such as the point (0,0) to measure the success of the algorithms used in the simulation, this may become a main issue if the optimal solutions are no long isolated points. We can expect that such types of functions can also make it more challenging for some algorithms to find their optimality.

\subsection{Function on a Hyperboloid}

We can easily extend a standard sphere function to a
sphere function with a hyperboloid (of revolution) constraint
\be f_5(\x)=\sum_{n=1}^D x_n^2, \ee
subject to
\be \sum_{n=1}^{D-1} \frac{x_n^2}{a^2} \ge \frac{x_D^2}{b^2} +1, \quad a, b \ge 1. \ee
It has a minimum $f_{\min}=a^2$ on a $(D-1)$-dimensional hyper-sphere
\be \sum_{n=1}^{D-1} x_n^2=a^2, \ee
which corresponds to infinitely many solutions.

In the case of 3D, we have
\be f(\x)=x^2+y^2+z^2, \ee
subject to \[ \frac{x^2}{a^2} +\frac{y^2}{a^2}-\frac{z^2}{b^2} \ge 1, \]
with $f_{\min}=a^2$ on $x^2+y^2=1$.

\subsection{Non-Smooth Multi-Layered Functions}

Though the above functions are more complicated, their objective values are still continuous in most cases. There is no jump or discontinuity in their landscapes.

In addition to make the search domains for independent decision variables irregular and/or to make the objectives with kinks, we can also design functions with discontinuous objective landscapes.
For example, we can make the objective values of a sphere function take only integer values, which can lead to a non-smooth multi-layered function
\be f(\x)=\Big\lfloor \sum_{n=1}^D x_n^2 \Big\rfloor, \ee
where $\lfloor x \rfloor$ is a floor function, which rounds $x$ to the nearest integer smaller than $x$.
That is, $k=\lfloor x \rfloor <=x$. For example, $\lfloor 2.3 \rfloor=2$ and $\lfloor 0.17 \rfloor =0$.

The optimal solution of this function is $f_{\min}=0$
inside a hyper-sphere
\be \sum_{n=1}^D x_n^2 =1. \ee
Any point inside this hyper-sphere is an optimal solution. Thus, this  function can have infinitely many optimal solutions within a hyper-volume, and all solutions inside this region have the same optimal objective value $f_{\min}=0$.

\begin{figure}
\begin{center}
\includegraphics[height=2.30in,width=4in]{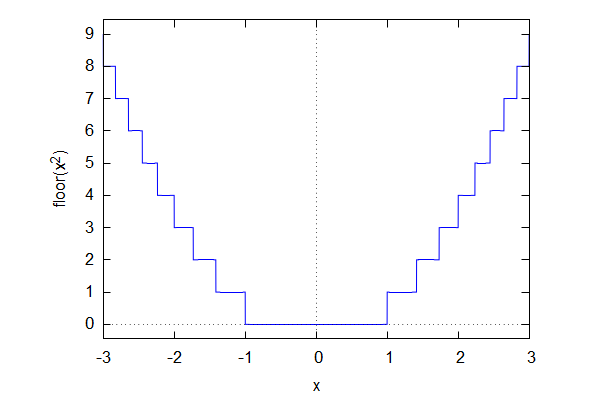}
\caption{A non-smooth function $f(x)=\lfloor x^2 \rfloor$ with an optimal region $-1 < x <+1$. \label{Fun-Fig-100}}
\end{center}
\end{figure}

On the one hand, it seems that this may be easier for algorithms to find the solution; however, many statistical measures such as mean solutions and standard deviations used for comparison would not make much sense for this function. Care should be taken when analyzing results. On the other hand, this integer-value objective function is no longer smooth, thus gradient-based methods would not work well.
For example, in the simplest one-dimensional case, the function is shown in Fig.~\ref{Fun-Fig-100}.

\begin{figure}
\begin{center}
\includegraphics[height=2in,width=3.5in]{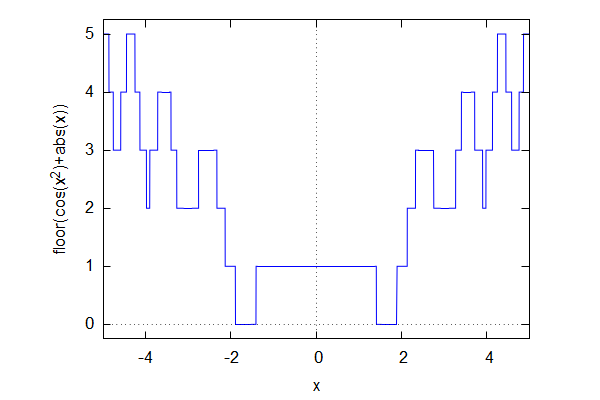}
\caption{Function $\lfloor |x|+\cos(x^2)\rfloor$ with two lowest flat regions that are both optimal. \label{Fig-floor-100} }
\end{center}
\end{figure}

To follow this same line of thinking, test benchmarks can have
multiple isolated regions as optimal regions with the same objective value. For example, we can design a function in the following form:
\be \min \;\; f_6(\x)=\left\lfloor \sum_{n=1}^D \Big[|x_n|+\cos(x_n^2)\Big] \right\rfloor. \ee
For $D=1$, we have
\be f(x)=\left\lfloor |x|+\cos(x^2)\right\rfloor, \ee
which has the global minimum value $f_{\min}=0$ with infinitely
many solutions in
two disconnected, flat regions: ($-1.89714$, $-1.41299$) and
($1.41299$, $1.89714$), which are shown in Fig.~\ref{Fig-floor-100}.

\begin{figure} \begin{center}
\includegraphics[height=2.5in, width=4in]{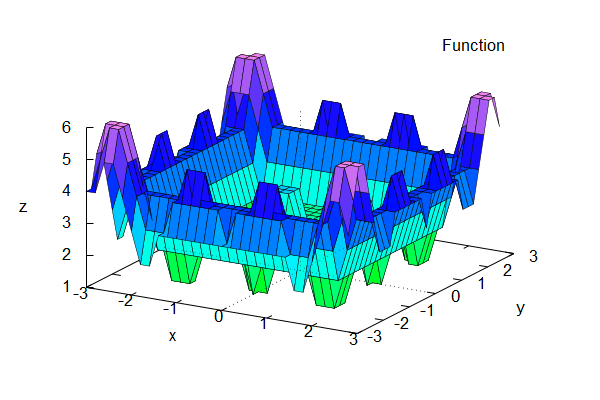}
\caption{Optimality occurs in multiple flat regions. \label{Fig-floor-200}}

\end{center} \end{figure}

In the case of $D=2$, we have
\be f(x,y)=\left\lfloor |x|+|y| +\cos(x^2) +\cos(y^2) \right\rfloor, \ee
which has 8 isolated flat regions with the same minimum $f_{\min}=0$,
as shown in Fig.~\ref{Fig-floor-200}.

\section{Parameter Estimation as Benchmarks}

For a vibration problem with a unit step input~\citep{YangMaths2017}, we have its mathematical equation as an ordinary differential
equation
\be \frac{\ddot{y}}{\omega^2} + 2 \zeta \frac{\dot{y}}{\omega} + y= u(t), \label{ODE-eq-100} \ee
where $\omega$ and $\zeta$ are the two parameters to be estimated.
Here, the unit step function is given
\be u(t)=\brak{0 & \quad \textrm{if } \;\; t<0, \\ 1 & \quad \textrm{if } \;\; t \ge 0. } \ee

Suppose for a given system, we have observed its actual response.
The relevant measurements are given Table~\ref{Table-data-100}.
\begin{table}
\begin{center}
\caption{Measured data for a vibration problem. \label{Table-data-100}}
\begin{tabular}{|r|rrrrrrrrrrrr|}
\hline
Time $t_i$ & 0 & & 1 & 2 & 3 & 4 & 5 & 6 & 7 & 8 & 9 & 10 \\
$y(t_i)$ & 0 & & 1.0706 & 1.3372 & 0.8277 & 0.9507 & 1.0848 & 0.9814 & 0.9769 & 1.0169 & 1.0012 & 0.9933 \\ \hline
\end{tabular}
\end{center}
\end{table}

In order to estimate the two unknown parameter values $\omega$ and $\zeta$, we can define the objective function as
\be f_7(\x)=\sum_{i=0}^{10} [y(t_i)-y_s(t_i)]^2, \ee
where $y(t_i)$ for $i=0,1,...,10$ are the observed values and
$y_s(t_i)$ are the values obtained by solving the differential equation~(\ref{ODE-eq-100}),
given a guessed set of values $\zeta$ and $\omega$.
Here, we have used $\x=(\zeta, \; \omega)$.

The true values are $\zeta=\frac{1}{4}$ and $\omega=2$. The aim of this benchmark is to solve the differential equation iteratively so as to find the best parameter values that minimize the objective or best-fit errors.

\section{Integrals as Benchmarks}

For almost all the benchmarks in the literature, integrals rarely appear in the problem formulations. Sometimes, the evaluation of an integral can be challenging, and thus the benchmarks involving integrals can make things more difficult.

From basic calculus, we know that it is very challenging to calculate the well-known integral
\be \int_0^{\infty} \frac{ \sin x}{x} dx =\frac{\pi}{2}. \ee
Suppose we want to maximize the integral
\be \max \;\; f_8(\beta, k)
=\int_0^{\infty} \frac{\sin (k x)}{ x e^{\beta x}} dx, \ee
where $\beta \ge 0$ is a real number and $k>0$ is an integer.

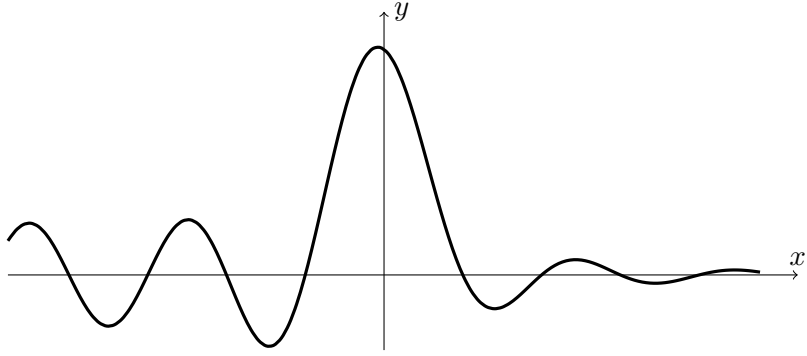
\begin{figure}
\begin{center}
\begin{tikzpicture}
\draw[very thick,domain=-5:5,samples=100] plot [smooth] (\x,{sin(deg(3*\x))/(\x*exp(0.25*\x)});
\draw[->] (-5,0)--(5.5,0) node [above] {$x$};
\draw[->] (0,-1)--(0,3.50) node [right] {$y$};
\end{tikzpicture}
\caption{Variation of $f(x)=\sin(k x)/(x \exp(\beta x))$ with $k=3$ and $\beta=0.25$. \label{Fig-sin-100}}
\end{center}
\end{figure}

Though theoretically we know that
\be \int_0^{\infty} \frac{\sin (k x)}{x e^{\beta x}} dx = \frac{\pi}{2}-\tan^{-1} \beta, \ee
which does not depend on $k$. The maximum value of this integral occurs at $\beta=0$ for any positive integer $k$.
In case of $\beta=0.5$ and $k=3$, the variation of the integrand is shown in Fig.~\ref{Fig-sin-100}.

Now let us use the maximization of this integral as a benchmark, which   requires to evaluate the integral with an infinite limit.
If without any prior knowledge of its true value of this integral, the effective evaluation of the objective function requires some sophisticated numerical integration. This requires a careful implementation.

\section{Benchmarks of Infinite Dimensions}

All the function benchmarks in the literature have reasonably
a well-defined solution set as the possible optimal solutions. Even in the high-dimensional space, such optimal solutions are just points or a region in the $D$-dimensional search space because the solutions are represented as a solution vector $\x$.

In some applications in science and engineering, the optimal solutions may not be represented by vectors. For example, in the
calculus of variations~\citep{YangMaths2017}, optimal solutions are curves and surfaces that cannot be represented by simple vectors.
For this branch of mathematics, we have rigorous theory using the Euler-Lagrange equation to solve such type of problems.

Now let us reformulate such problems and try to solve them using optimization algorithms without using the Euler-Lagrange equation. Since almost all optimization algorithms have been designed to  represent solutions in terms of vectors, this kind of problem from calculus of variations can be a major issue for standard optimization algorithms. In aerodynamics and shape optimization, shapes are represented by some parametric curves so as to simplify the representations. Even so, shape optimization can be quite challenging to implement.

To provide benchmarks for optimization algorithms from this perspective, let us use two examples to design benchmarks with
solutions as paths or curves. Even the actual space ($x,y$) is two-dimensional, the representation of solutions may require many points (or infinitely many points) to form a smooth curve. In fact,
we need a functional space to represent the potential curves properly. In this case, we can refer to this type of problem as benchmarks in infinite dimensions.

\subsection{Shortest-Path Problem}

In the two-dimensional space $(x,y)$, there is a path or curve $y(x)$ that minimize the integral
    \be \min \quad  Q=\int_0^1 \sqrt{1+\Big(\frac{dy}{dx}\Big)^2} dx. \ee
    From the plane geometry or calculus of variations~\citep{YangMaths2017}, we know that the solution is a straight line $y=x$ from the origin (0,0) to point (1,1).

    The challenge for this benchmark as an optimization problem is that the solution is not a single point but a segment of a curve (or a straight line in this case). In order to find this solution (corresponding to infinitely many points), some parametrization of an unknown curve is needed. Therefore, any standard algorithms for solving single-objective optimization problems have to be modified so that solution paths can be represented in effectively.

Mathematically, the objective functional can be rewritten as
    \be \textrm{minimize } \;\; f_9(y(x))=\int_0^1 \sqrt{1+\Big(\frac{dy}{dx}\Big)^2} dx, \ee
subject to $x \ge 0$ and  $y(x) \ge 0$.

\subsection{Shape Optimization}

The shape of a hanging rope under gravity takes the form that the potential energy is minimized, subject to a fixed length $L$.
The shape of a loose rope hinged between two fixed points $(-a,0)$ and $(a,0)$ can be obtained by minimizing the potential energy
\be \textrm{minimize } \;\; E_p=\rho g \int_{-a}^a y \sqrt{1+y'^2} \; dx, \ee
where $\rho$ is the density of the rope and $g$ is the acceleration due to gravity. Here, we use the notation $y'=dy/dx$ for a given smooth function $y(x)$. Without loss of generality, we can use $\rho g=1$, and thus we have
\be \textrm{minimize } \;\; f_{10}(y(x))=\int_{-a}^a y \sqrt{1+y'^2} \; dx, \ee
This minimization problem is subject to an equality
\be L=\int_{-a}^a \sqrt{1+y'^2} \; dx, \ee
with $L > 2a.$

\begin{figure}[h]
\begin{center}
\begin{tikzpicture}
\draw[very thick] (-4,0) .. controls (-2,-3.2) and (2, -3.2) .. (4,0);
\draw[thick] (-4,0) circle (0.05) node [above] {$(-a,0)$};
\draw[thick] (4,0) circle (0.05) node [above] {$(+a,0)$};

\end{tikzpicture}
\caption{Shape of a hanging rope with a fixed length. \label{Fig-shape-100}}
\end{center}
\end{figure}
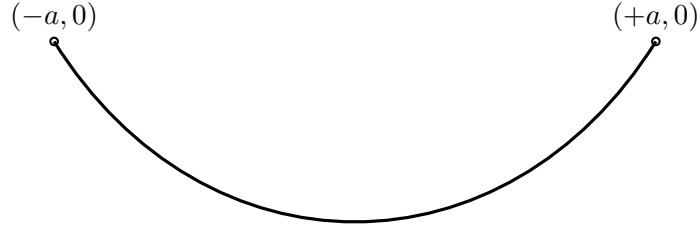

The solution can be obtained by solving the
Euler-Lagrange equation, which corresponds to
the shape or curve (see Fig.~\ref{Fig-shape-100})

\be y(x)=\cosh(x)-\cosh(a), \ee
with a length
\be L=\int_{-a}^a \sqrt{1+y'^2} \; dx =e^{a}-e^{-a}=2 \sinh(a). \ee
In case of $a=1$, we have $L=e-e^{-1} \approx 2.3504$.

The challenge of this benchmark is not only to represent the
solutions properly, but also to deal with the equality constraint correctly. In addition, the algorithm to be used must also be modified to accommodate such additional requirements.

\section{Conclusions}

There are many benchmarks for testing optimization algorithms; however, the benchmarks in the current literature tend to be smooth functions or problems with a finite number of optimal solutions. In addition, the search domains of existing benchmarks are usually regular with simple bounds or limits. To extend benchmarks to be more relevant to realistic problems, we have introduced ten new benchmarks by adding some noise, isolating the search domains and even making problems non-smooth with singularity and discontinuities.

These new benchmarks can be used to validate new algorithms and existing algorithms to see if they can cope with problems with non-smoothness and singularity well. We hope that this will inspire more research into benchmark problems and how to test new algorithms properly.


\begin{thebibliography}{100}

\bibitem{YangBook2014}
Yang, X.S.:
\newblock Nature-Inspired Optimization Algorithms.
\newblock Elsevier Insight, London (2014)

\bibitem{YangHe2019}
Yang, X.S., He, X.S.:
\newblock Mathematical Foundations of Nature-Inspired Algorithms.
\newblock Springer Briefs in Optimization. Springer, Cham, Switzerland (2019)

\bibitem{Yang2018Wiley}
Yang, X.S.:
\newblock Optimization Techniques and Applications with Examples.
\newblock John Wiley \& Sons, Hoboken, NJ, USA (2018)

\bibitem{Yang2020Rev}
Yang, X.S.:
\newblock Nature-inspired optimization algorithms: challenges and open
  problems.
\newblock Journal of Computational Science \textbf{Article 101104} (2020)

\bibitem{Yang2010Fun}
Yang, X.S.:
\newblock Firefly algorithm, stochastic test functions and design optimisation.
\newblock Int. J. Bio-Inspired Computation \textbf{2}(2) (2010)  78--84

\bibitem{Jamil2013}
Jamil, M., Yang, X.S.:
\newblock A literature survey of benchmark functions for global optimisation
  problems.
\newblock International Journal of Mathematical Modelling and Numerical
  Optimisation \textbf{4}(2) (2013)  150--194

\bibitem{Suganthan2005}
Suganthan, P., Hansen, N., Liang, J., Deb, K., Chen, Y., Auger, A., Tiwar, S.:
\newblock Problem definitions and evaluation criteria for cec 2005, special
  session on real-parameter optimization, technical report.
\newblock Technical report, Nanyang Technological University (NTU), Singapore
  (2005)

\bibitem{MazharFun}
Mazhar, A.A.:
\newblock Benchmark functions (web site).
\newblock Technical report, Victoria University of Wellington, New Zealand
  (2020)

\bibitem{HedarFun}
Hedar, A.:
\newblock Global optimization test problems (web site).
\newblock Technical report, University of Kyoto, Japan (2011)

\bibitem{Kumar2020}
Kumar, A., Wu, G., Ali, M.Z., Mallipeddi, R., Suganthan, P.N., Das, S.:
\newblock A test-suite of non-convex constrained optimization problems for the
  real-world and some baseline results.
\newblock Swarm and Evolutionary Computation \textbf{56}(Article 100693) (2020)

\bibitem{Zitzler2000}
Zitzler, E., Deb, K., Thiele, L.:
\newblock Comparison of multiobjective evolutionary algorithms: emperical
  results.
\newblock Evolutionary Computation \textbf{8}(2) (2000)  173--195

\bibitem{YangMaths2017}
Yang, X.S.:
\newblock Engineering Mathematics with Examples and Applications.
\newblock Academic Press, London (2017)

\end{thebibliography}
\end{document}